\documentclass[12pt]{article}
\usepackage{amsmath,amsfonts,latexsym,mathrsfs, amsthm,amssymb}
\makeatletter
  
  \@addtoreset{equation}{section}
\makeatother
\newtheorem{proposition}{Proposition}[section]
\newtheorem{theorem}{Theorem}[section]
\newtheorem{lemma}{Lemma}[section]
\newtheorem{remark}{Remark}[section]
\newtheorem{corollary}{Corollary}[section]
\newtheorem{definition}{Definition}[section]

\usepackage{color}
\usepackage{ulem}
\date{ }

\begin{document}

\title{Stochastic optimal transport with at most quadratic growth cost
}

\author{
Toshio Mikami\thanks{Partially supported by JSPS KAKENHI Grant Number 19K03548.
}\\
Department of Mathematics, Tsuda University
}

\maketitle

\begin{abstract}
We consider a class of stochastic optimal transport, SOT for short, with given two endpoint marginals
in the case where a cost function exhibits at most quadratic growth.
We first study the upper and lower estimates, the short--time asymptotics, the zero--noise limits, and the explosion rate as time goes to infinity of SOT.
We also show that the value function of SOT is equal to zero or infinity in the case where a cost function exhibits less than linear growth.
As a by--product, we characterize the finiteness of the value function of SOT by that of the Monge--Kantorovich problem with the same two endpoint marginals.
As an application, we show the existence of a continuous semimartingale, 
with given initial and terminal distributions,
of which the drift vector is $r$th integrable for $r\in [1,2)$.
We also consider the same problem for Schr\"odinger's problem where $r=2$.
This paper is a continuation of our previous work.
\end{abstract} 




\section{Introduction}
\label{intro}
For $d\ge 1, T>0$, let $\sigma :[0,T]\times \mathbb{R}^d\rightarrow
M(d,\mathbb{R})$ be a bounded Borel measurable $d\times d$--matrix function.
Let $\mathcal{A}_T$ denote the set of continuous semimartingales $\{ X(t)\}_{0\le t\le T}$ defined on a possibly different complete filtered probability space such that the following holds:
\begin{equation}\label{1.1Eq}
dX(t)= u_X (t)dt+\sigma(t,X(t))dB(t),\quad 0< t<T.
\end{equation}
Here  $\{u_X (t)\}_{0\le t\le T}$ and $\{B (t)\}_{0\le t\le T}$ are
a progressively measurable $\mathbb{R}^d$--valued stochastic process
and an $\mathbb{R}^d$--valued Brownian motion, respectively, defined on the same filtered probability space
 (see e.g., \cite{IW14}).
 If there exists $b:[0,T]\times \mathbb{R}^d\rightarrow\mathbb{R}^d$ such that
$u_X(\cdot)=b(\cdot,X(\cdot))$, 
then 
we write $b=b_{X}$.
 In this paper, the probability space under consideration is not fixed.
When it is not confusing, we use the same notation $P$ and $B$ for different probability and Brownian motion, respectively.
Let $\mathcal{P}(\mathbb{R}^d)$ denote the space of all Borel probability measures on $\mathbb{R}^d$ endowed with weak topology.
$$\mu_1(dx):=\mu(dx\times \mathbb{R}^d),\quad\mu_2(dx):=\mu(\mathbb{R}^d\times dx),\quad
\mu\in \mathcal{P}(\mathbb{R}^d\times \mathbb{R}^d).$$
$$\Pi(P,Q):=\{\mu\in \mathcal{P}(\mathbb{R}^d\times \mathbb{R}^d):
\mu_1=P,\mu_2=Q\},\quad P,Q\in \mathcal{P}(\mathbb{R}^d).$$
$$\mathcal{A}_T(P,Q):=\{\{ X(t)\}_{0\le t\le T}\in \mathcal{A}_T:
P^{(X(0), X(T))}\in \Pi(P,Q)\},$$ 
where $P^X$  denotes a probability distribution of a random variable $X$.
Let $L:[0,T]\times \mathbb{R}^d \times \mathbb{R}^d\longrightarrow [0,\infty )$ be  Borel measurable.

In this paper, we consider the upper and lower estimates and the short and long--time asymptotics of 
the following Stochastic optimal transport, SOT for short 
 (see \cite{M08, M09, M21, M2021} for SOT and related topics).
 
\begin{definition}[Stochastic optimal transport\label{def1}]

\noindent
For $P,Q\in \mathcal{P}(\mathbb{R}^d ), T>0$,
\begin{align}
V(T,P,Q)
&:=\inf \left\{E\biggl[\int_0^T L(t,X(t);u_X (t))dt \biggr]:X\in \mathcal{A}_T(P,Q)\right\},\\
V^{fc} (T,P,Q)
&:=\inf \left\{E\biggl[\int_0^T L(t,X (t);b_X(t,X(t)))dt \biggr]:\right.\nonumber \\
&\qquad \left.X\in \mathcal{A}_T(P,Q),  u_X(\cdot)=b_X(\cdot,X(\cdot))\right\}.\label{12}
\end{align}
When $L=|u|^r, r>0$, we write $V=V_r$.
If the set over which the infimum is taken is empty, we set the infimum infinite.
\end{definition}

\begin{remark}\label{Remark1.1}
(i) ``fc'' in (\ref{12}) means ``feedback control'' (see e.g., \cite{FS06}).
(ii) $V^{fc}\ge V$.
Equality holds if
$u\mapsto L(t,x;u)$ is convex and is of at least linear growth,  uniformly in $(t,x) \in [0,T]\times \mathbb{R}^d$, in which case we can assume, without loss of generality, that $u_X(t)$ is $(\mathcal{F}^{X}_t:=\sigma [X(s), 0\le s\le t])$--measurable (see \cite {M21, M2021}).
(iii) 
If $u\mapsto L(t,x;u)$ is convex for $(t,x) \in [0,T]\times \mathbb{R}^d$, is of more than linear growth uniformly in (t,x), 
$(t,x,u)\mapsto L(t,x;u)$ is lower semicontinuous, and $V(T, P,Q)$ is finite,
then there exist a minimizer of $V (T,P,Q)$.
In particular, there exists a continuous semimartingale $\{X(t)\}_{0\le t\le T}$  and a function $b_X$ such that the following holds 
(see \cite {M21} and the references therein):
\begin{align}
dX(t)&=b_X(t,X(t))dt+\sigma (t,X(t))dB(t),\quad 0<t<T,\label{1.9}\\
P^{(X(0),X (T))}&\in \Pi(P,Q).\label{1.10}
\end{align}
In addition, if $|u|^r/(1+L(t,x;u))$ is bounded, 
then the following holds:
\begin{equation}
E\biggl[\int_0^T |b_X(t,X(t))|^r dt\biggr]<\infty.\label{1.11}
\end{equation}
\end{remark}

$a:=\sigma \sigma^*$, where  $\sigma^*$ denotes the transpose of $\sigma$.
\begin{align*}
||A||&:=\left(\sum_{i=1}^dA_{ij}^2\right)^{1/2},\quad A=(A_{ij})_{i,j=1}^d\in M(d,\mathbb{R}),\\
||A||_{\infty,T}&:=\sup\left\{||A(t,x)||: (t,x)\in [0,T]\times \mathbb{R}^d\right\},\\
&\qquad\qquad A:[0,T]\times \mathbb{R}^d\rightarrow M(d,\mathbb{R}),\\
|f|_{\infty,T}&:=\sup\{|f(t,x)|: (t,x)\in [0,T]\times \mathbb{R}^d\},\quad %
f:[0,T]\times \mathbb{R}^d\rightarrow \mathbb{R}^d,\\
\mathcal{P}_r (\mathbb{R}^d )&:=\left\{P\in \mathcal{P} (\mathbb{R}^d ):
\int_{\mathbb{R}^d}|x|^rP(dx)<\infty\right\},\quad r>0.
\end{align*}
Consider $\{x+s(y-x)/t\}_{0\le s\le t}$ instead of $\{x+s(y-x)\}_{0\le s\le 1}$ in  \cite{M15}.
Then the following holds.

\begin{theorem}[Lemma 2.1 and Corollary 2.4 in \cite{M15}]\label{thm1.2}
Suppose that $\sigma (\cdot)=(\sigma_{ij}(\cdot))_{i,j=1}^d$ is uniformly nondegenerate on $[0,T]\times \mathbb{R}^d$ and 
$\sigma_{ij}\in C^{1,2}_b ([0,T]\times \mathbb{R}^d)$, $i,j=1,\cdots ,d$.
Then for any $r>1$ and any $P,Q\in \mathcal{P}_r(\mathbb{R}^d)$ such that 
$p(x):=P(dx)/dx$ and $q(x):=Q(dx)/dx$ exist and are absolutely continuous in $x$,
the following holds: for $t\in (0,T]$,
\begin{align}\label{110}
&t^{r-1}V_r(t,P,Q)\nonumber\\
&\le
2^{r-2}\left(2^r\int_{\mathbb{R}^d\times \mathbb{R}^d}\left|y-x\right|^rp(x)q(y)dxdy
+t^r|(D_x^*a)^*|_{\infty,t}^r\right)\nonumber\\
&\qquad +2^{r-1}t^r||a||_{\infty,t}^r\sum_{f=p,q}
\int_{\mathbb{R}^d}|D_x\log f(x)|^rf(x)dx.
\end{align}
In particular, if 
\begin{equation}\label{111}
\sum_{f=p,q}\int_{\mathbb{R}^d}|D_x\log f(x)|^r f(x)dx<\infty,
\end{equation}
then $\{t^{r-1}V_r(t,P,Q)\}_{0<t\le T}$ is bounded and there exist $b_X$ and $X$ such that (\ref{1.9})--(\ref{1.11}) hold. 
Besides,
\begin{equation}
\limsup_{t\to 0} t^{r-1}V_r(t,P,Q)\le
4^{r-1}\int_{\mathbb{R}^d\times \mathbb{R}^d}|y-x|^rp(x)q(y)dxdy.
\end{equation}
\end{theorem}

In this paper, for $r\in [1,2)$, we generalize Theorem \ref{thm1.2} without   
an assumption (\ref{111}).
As a by--product, we show the existence of $X\in \mathcal{A}_T$ for which  (\ref{1.9})--(\ref{1.11}) hold
even when the h--path process does not exist (see Proposition \ref{Thm2.1} and Corollary \ref{corollary2.1}  in section \ref{sec:2.1}).
For $r\ge 1$, we also give  the lower bound of $V_r(t,P,Q)$ (see Proposition \ref{pp2.2.0} in section \ref{sec:2.1}).
In particular, we show that $\{t^{r-1}V_r(t,P,Q)\}_{0<t\le T}$ is bounded if and only if 
\begin{equation}\label{19}
{\bf T}_r (P,Q):=\inf \left\{\int_{\mathbb{R}^d \times \mathbb{R}^d} |y-x|^r\mu(dxdy):
\mu\in \Pi (P,Q)\right\}<\infty
\end{equation}
(see e.g., \cite{1-RR, Vi} for optimal transport).
Besides, we show that the following holds: if  ${\bf T}_r (P,Q)$ is finite,
\begin{align*}
\lim_{t\to 0}t^{r-1}V_r(t,P,Q)&={\bf T}_r (P,Q),\\
\limsup_{t\to 0}t^{r/2-1}V_r(t,P,P)&<\infty, \\
\lim_{t\to 0}\{V_r(t,P,Q)-t^{1-r}{\bf T}_r (P,Q)\}&=0,\quad1\le r<\frac{3}{2},\\
\lim_{t\to 0}\{V_r(t, P,Q)^{1/r}-(t^{1-r}{\bf T}_r (P,Q))^{1/r}\}&=0,\quad1\le r<2
\end{align*}
 (see Corollary \ref{co2.3} in section \ref{sec:2.1} and also \cite{1-ADPZ,1-DLR,L2,M04,PAL}
 for the related problem in Schr\"odinger's problem).
 $V(t,P,Q)$ plays a role of a fundamental solution for the Nisio semigroup.
Our result gives a version of Aronson--type inequality for $\exp (-V_r(t,P,Q))$ with $|x-y|^2/t$ replaced by ${\bf T}_r (P,Q)/t^{r-1}$, provided $r\in [1, 2)$ (see \cite{A67, Sheu91}).
We also study the zero--noise limit, and the explosion rate as time goes to infinity of $V_r(t,P,Q)$ (see Corollary \ref{co2.2} and Theorem \ref{Theorem2.1} in section \ref{sec:2.1}).
 For $r\in (0,1)$, we also show that $V_r(t,P,Q)=0$ (see Proposition \ref{pp2.1} in section \ref{sec:2.1}).

Let $\xi(t,x):[0,T]\times \mathbb{R}^d\longrightarrow \mathbb{R}^d$ be Borel measurable.
The following is a typical example of SOT (see \cite{B32,Csiszar,F88,L2,M08,M09,M2021,S1,S2,Z86} and the references therein).

\begin{definition}[Schr\"odinger's problem\label{def2}]
We denote $V$ by $V^S$ and call it  Schr\"odinger's problem
if the following holds:
\begin{equation}\label{1.2}
L(t,x;u)=\frac{1}{2}|\sigma^{-1}(t,x)(u -\xi (t,x))|^2,\quad (t,x,u)\in [0,T]\times \mathbb{R}^d \times \mathbb{R}^d.
\end{equation}
\end{definition}

We introduce the following.

\noindent
(A) 
$a(t,x)$, $(t,x)\in [0,T]\times \mathbb{R}^d$, is uniformly nondegenerate,
bounded, once continuously differentiable, and uniformly H\"older continuous.
$D_x a(t,x)$ is bounded and the first derivatives of $a(t,x)$ are uniformly H\"older continuous
in $x$ uniformly in $t\in [0,T]$.
$\xi(t,x)$, $(t,x)\in [0,T]\times \mathbb{R}^d$ is bounded, continuous, and uniformly H\"older continuous
in $x$ uniformly in $t\in [0,T]$.

The following is known.

\begin{theorem}[see \cite{J75}]\label{thm11}
Suppose that (A) holds. Then the following SDE has a unique weak solution 
with a positive continuous transition probability density $p(s,x;t,y), 0\le s<t\le T, x,y\in \mathbb{R}^d$:
 for $P\in \mathcal{P}(\mathbb{R}^d )$,
\begin{align}\label{1.4}
dX(t)&=\xi(t,X(t))dt+\sigma (t,X(t))dB(t),\quad 0<t<T,\\
P^{X(0)}&=P.\nonumber
\end{align}
Besides, for $Q\in \mathcal{P}(\mathbb{R}^d )$  such that $Q(dx)\ll dx$,
there exists a unique $\sigma$--finite measure $\nu_1(dx)\nu_2(dy)$ such that  (\ref{1.9})--(\ref{1.10}) with 
\begin{equation}\label{115}
b_X(t,x)=\xi(t,x)+a(t,x)D_x\log \int_{\mathbb{R}^d} p(t,x;T,y)\nu_2(dy)
\end{equation}
has a unique weak solution $X_o$ and that the following holds:
\begin{equation}\label{113}
P^{(X_o(0),X_o(T))}(dxdy)=\nu_1(dx)p(0,x ;T,y)\nu_2(dy).
\end{equation}
\end{theorem}

$X_o$ in Theorem \ref{thm11} is called the h--path process or the Schr\"odinger process.
If $V^S(T,P,Q)$ is finite, then, under (A), 
\begin{align}\label{1.15}
V^S(T,P,Q)&=v^S(T,P,Q)\nonumber\\
&:=\inf\{H(\mu(dxdy)|P(dx)p(0,x;T,y)dy)|\mu\in \Pi (P,Q)\}.
\end{align}
Here 
$$H(\mu|\nu):=
\begin{cases}
\displaystyle\int_{\mathbb{R}^d\times \mathbb{R}^d} \left(\log\frac{\mu (dx)}{\nu(dx)} \right)\mu (dx),&\mu\ll\nu,\\
\infty,&\hbox{otherwise,}
\end{cases}
\quad \mu,\nu\in \mathcal{P}(\mathbb{R}^d\times \mathbb{R}^d).
$$
$X_o$ and $P^{(X_o(0),X_o(T))}$ are 
unique minimizers of $V^S(T,P,Q)$ and of $v^S(T,P,Q)$, respectively
(see \cite{Csiszar,S1,S2,Z86,Z86-2,1-Zambrini3} and also \cite{D91,F88,L2,M08,M2021,MT06,RT93} and the references therein).
We  write $\mu_T(P,Q):=P^{(X_o(0),X_o(T))}$ when we consider $v^S(T,P,Q)$
and also $\mu_T(P,Q)=\mu_T$ for simplicity, when it is not confusing.
The continuities of $\mu_T(dxdy)$ and $\nu_1(dx)\nu_2(dy)$ in $(P,Q, p)$ in strong and weak topologies are given in \cite{M19} and \cite{mikami2021}, respectively.
The displacement semi--concavity and the Lipschitz continuity in the 2-Wasserstein distance ${\bf T}_2^{1/2}$
of $P\mapsto V^S(T,P,Q)$ are given in \cite{M21} (see also \cite{M2021}).
The weak continuity of $(P,Q, p)\mapsto \mu_T(dxdy)$ plays a crucial role in this paper.
\begin{equation}\label{16}
\mathcal{S}(P):=
\begin{cases}
\displaystyle\int_{\mathbb{R}^d}p(x)\log p(x)dx, &P(dx)=p(x) dx,\\
+\infty,&\hbox{otherwise,}
\end{cases}
\quad P\in \mathcal{P} (\mathbb{R}^d ).
\end{equation}
The following is known (see e.g., Theorem 8 in \cite{M21} and Lemma \ref{lm31} in section \ref{sec:3}).

\begin{theorem}\label{thm1.1}
Suppose that (A) holds. 
Then for any $P,Q\in  \mathcal{P}_{2} (\mathbb{R}^d )$, 
$\{tV^S(t,P,Q)\}_{0<t\le T}$ is bounded if and only if 
$\mathcal{S}(Q)$ is finite.
\end{theorem}

\begin{remark}\label{rk1.2}
(i) (\ref{111}) with $r=2$ implies that $\mathcal{S}(Q)$ is finite.
Indeed, by the log-Sobolev inequality for Gaussian measures (see e.g., \cite {Vi}),
the following holds: for $p(x)dx\in \mathcal{P}_2(\mathbb{R}^d)$ such that
$x\mapsto p(x)$ is absolutely continuous,
$$
\int_{\mathbb{R}^d}p(x)\log p(x)dx
\le 
\int_{\mathbb{R}^d}p(x)\log g(x)dx+\frac{1}{2}\int_{\mathbb{R}^d}|D_x\log p(x)+x|^2 p(x)dx.
$$
Here 
$$g(x):=\frac{1}{\sqrt{2\pi }^d}\exp \left(-\frac{|x|^2}{2}\right),\quad 
x\in  \mathbb{R}^d.$$
(ii) For $p(x)dx\in \mathcal{P}_2(\mathbb{R}^d)$,
$\log p(x)\in L^1(p(x)dx)$ if and only if $H(p(x)dx|g(x)dx)$ is finite.
\end{remark}
In this paper, we also study the upper and lower estimates
on a finite time interval and the short--time asymptotics of $V^S(t,P,Q)$.
In the case where there exists an invariant probability  density $m(y)$ of the solution to the SDE (\ref{1.4}) 
and where the transition probability density converges to $m(y)$, $t\to\infty$, locally uniformly,
we discuss the upper  bound of $V^S(t,P,Q)$ on an infinite time interval and the long--time asymptotics of $V^S(t,P,Q)$ and  $\mu_t (P,Q)$
(see Proposition \ref{pp2.2}, Theorem \ref{pp2.3} and Corollary \ref{co2.5} in section \ref{sec:2.2}).
The long--time asymptotics of $V^S(t,P,Q)$ and $\mu_t (P,Q)$ were studied in \cite{Con} 
when $P$ and $Q$ have compact supports (see Remark \ref{Remark2.3} for more details).

In section  \ref{sec:2}, we state our results and prove them in sections \ref{sec:3}--\ref{sec:4}.

\section{Main result} \label{sec:2}

In this section, we state our results.

We introduce our assumptions.

\noindent
(A1) 
(i) $\sigma\in C_b([0,T]\times \mathbb{R}^d;M(d,\mathbb{R}))$.
$x\mapsto \sigma (t,x)$ is uniformly Lipschitz continuous, uniformly in $t\in [0,T]$.
(ii) $\sigma(t,x)=\sigma (t,0)=:\sigma (t),  (t,x)\in [0,T]\times \mathbb{R}^d$. 

\noindent
(A1)' $0<\lambda_m\le \lambda_M<\infty$, where
\begin{align*}
\lambda_m&:=\inf\{\hbox{Trace} (a (t,x)):(t,x)\in [0,\infty)\times \mathbb{R}^d\},\\
\lambda_M&:=\sup\{\hbox{Trace}(a (t,x)):(t,x)\in [0,\infty)\times \mathbb{R}^d\}.
\end{align*}

\noindent
(A2) 
(i) There exist $r\in [1,2)$ and $C_{r,t}, C_{r,T}'>0 $ such that 
$\{C_{r,t}\}_{0\le t\le T}$ is bounded and that the following holds:
\begin{equation}\label{1.19}
L(t,x;u)\le C_{r,t}|u|^r+C_{r,T}',\quad (t,x,u)\in [0,T]\times \mathbb{R}^d\times \mathbb{R}^d.
\end{equation}

\noindent
(ii) There exist $r\ge 1$ and $c_{r,t}, c_{r,T}'>0$ such that the following holds:
\begin{equation}
L(t,x;u)\ge c_{r,t}|u|^r-c_{r,T}',\quad (t,x,u)\in [0,T]\times \mathbb{R}^d\times \mathbb{R}^d.
\end{equation}

\noindent
(A3)
$L(t,x;0)=0, (t,x)\in [0,T]\times \mathbb{R}^d$.
There exists $r\in (0,1)$  and $C_{r,T}, C_{r,T}'>0$ such that (\ref{1.19}) 
with $C_{r,t}$ replaced by $C_{r,T}$ holds.

\begin{remark}
If $\sigma$ is bounded, then, for $r\ge 1$,
$V_r(T, P,Q)$ is equal to $T^{1-r}V_r(1, P,Q)$ with $a(t,\cdot), 0\le t\le T$ replaced by 
$Ta(Tt,x), 0\le t\le 1$ (see Proposition 1 in \cite{M21}).
\end{remark}

When $\sigma=\sigma (x), \xi=\xi (x)$ in (A) in section \ref{intro}, we denote it by (A0).

\noindent
(A0) $a :\mathbb{R}^d\rightarrow M(d,\mathbb{R})$ is uniformly nondegenerate,
 $a\in C_b^1(\mathbb{R}^d)$, and $Da$ is uniformly H\"older continuous.
$\xi:\mathbb{R}^d\longrightarrow \mathbb{R}^d$ is bounded and uniformly H\"older continuous.
 
$$r(\ell):=\frac{1}{2}\left\{\sup_{x\in \mathbb{R}^d}\{\hbox{Trace }a(x)\}+(\ell-2)
\sup_{x\in \mathbb{R}^d, x\ne 0} \left\langle a(x)\frac{x}{|x|},\frac{x}{|x|}\right\rangle\right\}.$$
(A4) 
There exist $M>0, \ell>2, r>r(\ell)$ such that the following holds:
$$\langle \xi(x),x\rangle\le -r,\quad|x|\ge M.$$

\begin{remark}\label{remark21}
Let $p(t,x,y)$ denote the transition probability density of the solution to the SDE (\ref{1.4}) when (A0) holds.
Under (A0) and (A4), the solution to (\ref{1.4}) has a unique invariant probability  density
$m$ (see \cite{V-2}).
$m\in C_b (\mathbb{R}^d;(0,\infty))$ since $p(1,\cdot, \cdot)\in C_b (\mathbb{R}^d\times \mathbb{R}^d;(0,\infty))$
(see Theorem \ref{thm11} and Lemma \ref{lm31}), and
$$m(x)=\int_{\mathbb{R}^d}m(y)p(1,y,x)dy.$$
\end{remark}

We first consider the SOT with less than quadratic growth cost in section \ref {sec:2.1} and then Schr\"odinger's problem in section \ref{sec:2.2}.

\subsection{SOT with less than quadratic growth cost}\label{sec:2.1}
 
The following generalizes  (\ref{110}) when $r\in (0,2)$.
\begin{proposition}\label{Thm2.1}
Suppose that (A1,i) holds.
Then for any $r\in (0,2)$, $P,Q\in  \mathcal{P} (\mathbb{R}^d )$ and $t\in (0,T]$,
\begin{equation}\label{2.2.0527}
t^{r-1}V_r( t,P,Q)\le {\bf T}_r (P,Q)+
\frac{2||\sigma||_{\infty,t}^r }{2-r}t^{r/2}.
\end{equation}
\end{proposition}

The following gives  a lower bound of $V_r(t, P,Q)$.

\begin{proposition}\label{pp2.2.0}
Suppose that (A1,i) holds.
Then for any $r\ge 1$, there exists $C_r>0$ such that the following holds:  for any $P,Q\in  \mathcal{P}(\mathbb{R}^d )$ and any $t\in (0,T]$ and $\varepsilon\in (0,1)$,
\begin{align}
(t^{r-1}V_r(t,P,Q))^{1/r}
&\ge{\bf T}_r(P,Q)^{1/r}-(C_r)^{1/r}||\sigma ||_{\infty,t}t^{1/2},\label{2.7.0}\\
t^{r-1}V_r(t,P,Q)&\ge (1-\varepsilon)^{r-1}{\bf T}_r(P,Q)-\varepsilon^{1-r}(1-\varepsilon)^{r-1}C_r||\sigma ||_{\infty,t}^rt^{r/2},\label{2.6}\\
(t^{r-1}V_r(t,P,Q))^{1/r}
&\ge {\bf T}_r(P^{X^0(t)},Q)^{1/r}-(C_r)^{1/r}t^{1/2}\label{2.5}\\
&\quad \times\sup\{||\sigma (s,x)-\sigma (s,y)||:0\le s\le t, x,y\in\mathbb{R}^d\},\nonumber\\
t^{r-1}V_r(t,P,Q)
&\ge  (1-\varepsilon)^{r-1}{\bf T}_r(P^{X^0(t)},Q)-\varepsilon^{1-r}(1-\varepsilon)^{r-1}C_rt^{r/2}\label{2.4.0}\\
&\qquad\times\sup\{||\sigma (s,x)-\sigma (s,y)||^r:0\le s\le t, x,y\in\mathbb{R}^d\}.\nonumber
\end{align}
Here $X^0$ denotes $X$ with $u_X=0$, $P^{X^0(0)}=P$
and $C_r$ can be replaced by $1$ for $r\in [1,2]$.
In particular, if, in addition, (A1,ii) holds, then
\begin{equation}\label{2.6.0}
t^{r-1}V_r(t,P,Q)\ge {\bf T}_r(P^{X^0(t)},Q).
\end{equation}
\end{proposition}

\begin{remark}
Suppose that the following holds:
\begin{equation}\label{2.9.0906}
{\bf T}_r(P,Q)\le t^{r-1}V_r(t,P,Q)+ C_r||\sigma ||_{\infty,t}^rt^{r/2}.
\end{equation}
Then (\ref{2.7.0}) and then (\ref{2.6}) hold (see (\ref{4.9.m1}) and (\ref{4.6.0905})).
We study if this inequality holds in future since (\ref{2.2.0527}) and (\ref{2.9.0906}) imply that 
(\ref{2.8.527}) holds for all $r\in [1,2)$. 
\end{remark}

From Propositions \ref{Thm2.1}--\ref{pp2.2.0}, the following holds.

\begin{corollary}\label{co2.3}
Suppose that (A1,i)  holds.
Then for any $r\in [1,2)$ and $P,Q\in  \mathcal{P} (\mathbb{R}^d )$, 
$\{t^{r-1}V_r(t, P,Q)\}_{0< t\le T}$ is bounded if and only if ${\bf T}_r (P,Q)$ is finite. 
Besides, if ${\bf T}_r (P,Q)$ is finite, then the following holds.
\begin{equation}\label{2.7}
\limsup_{t\to 0}\{t^{-1/2}|t^{r-1}V_r(t, P,Q)-{\bf T}_r (P,Q)| \}
\le 2r||\sigma ||_{\infty,0+}{\bf T}_r(P,Q)^{1-1/r},
\end{equation}
where $0^0:=1$.
\begin{equation}\label{2.8}
\limsup_{t\to 0}\{t^{r/2-1}V_r(t,  P,P)\}\le\frac{2||\sigma||_{\infty,0+}^r }{2-r}. 
\end{equation}
In particular, 
\begin{equation}\label{2.8.527}
\lim_{t\to 0}\{V_r(t,  P,Q)-t^{1-r}{\bf T}_r (P,Q)\}=0,\quad r\in [1,3/2).
\end{equation}
\begin{equation}\label{2.11}
\lim_{t\to 0}\{t^{-1/2}|(t^{r-1}V_r(t, P,Q))^{1/r}-{\bf T}_r (P,Q)^{1/r}| \}
\le \left(\frac{2}{2-r}\right)^{1/r} 
||\sigma ||_{\infty,0+}.
\end{equation}
In particular, 
\begin{equation}\label{2.12}
\lim_{t\to 0}\{V_r(t, P,Q)^{1/r}-(t^{1-r}{\bf T}_r (P,Q))^{1/r}\}=0,\quad r\in [1,2).
\end{equation}
\end{corollary}
\begin{remark}
(\ref{2.8})  implies  (\ref{2.7}) when $P=Q$ since $r\ge 1$.
(\ref{2.11})--(\ref{2.12}) do not improve (\ref{2.8.527}).
Indeed, by the mean value theorem,  there exists $\theta\in (0,1)$ such that
the following holds (see (\ref{2.12})):
\begin{align*}
&|V_r(t,  P,Q)-t^{1-r}{\bf T}_r (P,Q)|\\
&=r\{\theta V_r(t,  P,Q)^{1/r}+(1-\theta )(t^{1-r}{\bf T}_r (P,Q))^{1/r}\}^{r-1}\\
&\qquad \times |V_r(t,  P,Q)^{1/r}-(t^{1-r}{\bf T}_r (P,Q))^{1/r}|\\
&\approx 
t^{3/2-r} {\bf T}_r (P,Q)^{(r-1)/r}\times t^{-1/2} |(t^{r-1}V_r(t,  P,Q))^{1/r}-{\bf T}_r (P,Q)^{1/r}|.
\end{align*}
\end{remark}

We easily obtain the following from Corollary \ref{co2.3} above and Proposition 1 and its proof in \cite{M21}, where  D. Trevisan's superposition principle played a crucial role.
We omit the proof.

\begin{corollary}\label{corollary2.1}
Suppose that (A1,i) holds.
Then for any $r\in [1,2)$ and for any $P,Q\in  \mathcal{P} (\mathbb{R}^d )$ such that ${\bf T}_r (P,Q)$ is finite, there exist $X$ and  $b_X$  such that (\ref{1.9})--(\ref{1.11}) hold.
\end{corollary}

\begin{remark}
If ${\bf T}_r (P,Q)$ is finite for some $r\ge 2$, then ${\bf T}_r (P,Q)$ is finite for $r\in (1,2)$.
In particular, our result is meaningful when $\mathcal{S} (Q)$ is not finite
(see Theorem \ref{thm1.1} in section \ref{intro}).
\end{remark}

Under (A2), we immediately obtain an upper and lower bounds for $V(t,P, Q)$ from Propositions \ref{Thm2.1}--\ref{pp2.2.0} and Corollary \ref{co2.3} (see Remark \ref{Remark1.1}, (i)--(ii)).

\begin{corollary}
(i) Suppose  that (A1,i) and (A2,i) hold.
Then for any $P,Q\in  \mathcal{P} (\mathbb{R}^d )$ such that ${\bf T}_r (P,Q)$ is finite,
$\{t^{r-1}V^{fc}(t, P,Q)\}_{0< t\le T}$ is bounded.
Besides, we obtain the upper bound of $V^{fc}(t, P,Q)$ via that of  $V_r(t, P,Q)$.
In particular,
\begin{align} 
\limsup_{t\to 0}\{t^{r-1}V^{fc}(t,P,Q)\}&\le  \left(\limsup_{t\to 0}C_{r,t}\right) {\bf T}_r (P,Q),\\
\limsup_{t\to 0}\{t^{r/2-1}V^{fc}(t,P,P)\}&\le \left(\limsup_{t\to 0}C_{r,t}\right) \frac{2||\sigma||_{\infty,0+}^r}{2-r}. 
\end{align} 
(ii) Suppose that (A1,i) and (A2,ii) hold.
Then for any $P,Q\in  \mathcal{P} (\mathbb{R}^d )$,
we obtain the lower bound of $V(t, P,Q)$ via that of  $V_r(t, P,Q)$.
In particular,
\begin{equation} 
\liminf_{t\to 0}\{t^{r-1}V(t,P,Q)\}\ge \left(\liminf_{t\to 0}c_{r,t}\right) {\bf T}_r (P,Q).
\end{equation} 

\end{corollary}

For $\varepsilon>0$, let $V^\varepsilon_r$ denote $V_r$ with $\sigma$ replaced by
$\sqrt\varepsilon \sigma$.
The zero--noise limit of $V^\varepsilon_r$ can be obtained exactly in the same way as 
Corollary \ref{co2.3} and we omit the proof.
Indeed, for $V^\varepsilon_r$, on the r.h.s. of  in (\ref{2.2.0527}) and (\ref{2.6}),
$||\sigma||_{\infty,t}^r$ and $t^{r/2}$ can be replaced by
$t^{1/2}||\sigma||_{\infty,t}^r$ and $\varepsilon^{r/2}$, respectively
since
$$||\sqrt\varepsilon \sigma||_{\infty,t}^r t^{r/2}=\{t^{1/2}||\sigma||_{\infty,t}\}^r \varepsilon^{r/2}.$$

\begin{corollary}\label{co2.2}
Suppose that (A1,i)  holds.
Then for any $r\in [1,2)$ and $P,Q\in  \mathcal{P} (\mathbb{R}^d )$ for which ${\bf T}_r (P,Q)$ is finite, the following holds.
\begin{equation}
\limsup_{\varepsilon\to 0}\{\varepsilon^{-1/2}|V^\varepsilon_r(T, P,Q)-T^{1-r}{\bf T}_r (P,Q)| \}
\le 2r||\sigma ||_{\infty,T}T^{3/2-r}{\bf T}_r(P,Q)^{1-1/r},
\end{equation}
\begin{equation}
\limsup_{\varepsilon\to 0}\{\varepsilon ^{-r/2}V^\varepsilon_r(T,  P,P)\}\le\frac{2||\sigma||_{\infty,T}^r T^{1-r/2}}{2-r}. 
\end{equation}
\end{corollary}

The following gives the explosion rate of $V_r(t,P,Q)$ as $t\to\infty$.

\begin{theorem}\label{Theorem2.1}
Suppose that (A1)' holds and that $r\in [1,2)$.
Then there exist $C(r,\lambda_m,\lambda_M), C'(r,\lambda_M)>0$ such that for any $P,Q\in \mathbb{P}_r(\mathbb{R}^d)$ and $t>0$, the following holds:
\begin{align}
&V_r(t,P,Q)\nonumber\\
&\ge t^{1-r/2}
\left\{C(r,\lambda_m,\lambda_M)-t^{-1/2}C'(r,\lambda_M)
\left(\int_{\mathbb{R}^d}|x|^r(P(dx)+Q(dx))\right)^{1/r}\right\}\nonumber\\
&\to \infty,\quad t\to\infty.
\end{align}
\end{theorem}

For $r\in (0,1)$, the following holds (see \cite{MY22} for a similar problem in optimal transport).

\begin{proposition}\label{pp2.1}
Suppose that (A1,i) and (A3) hold.
Then for any $P,Q\in  \mathcal{P} (\mathbb{R}^d )$ for which $T_r (P,Q)$ is finite,
$V(T,P,Q)=0$.
Suppose, in addition, that (A1,ii) holds.
Then for any $P,Q\in  \mathcal{P} (\mathbb{R}^d )$  for which $V_r (T, P,Q)$ is finite,
$V(T,P,Q)=0$.
\end{proposition}

\begin{remark}
From Proposition \ref{Thm2.1}, under (A1,i), for $r\in (0,1)$, if $T_r (P,Q)$ is finite,
then $V_r (T, P,Q)$ is finite.
\end{remark}

\subsection{Schr\"odinger's problem}\label{sec:2.2}

In this section, we consider Schr\"odinger's problem.

 Let $\lambda_{\infty,t}$ denote the supremum of maximal eigenvalue of $a(s,x)$ on $[0,t]\times \mathbb{R}^d$.
The following implies the continuity and gives upper and lower bounds of $V^S(\cdot,P,Q)$ on a finite time interval.
It also gives the explosion rate of $V^S(t,P, Q)$ as $t\to 0$ when $P\ne Q$
(see \cite{L2, M04,M2021}).
We give the proof  for the sake of completeness.
(\ref{2.4}) is a rough estimate compared to Proposition \ref{Thm2.1} where $0<r<2$.
Recall that  Proposition \ref{pp2.2.0} holds for all $r\ge 1$.

\begin{proposition}\label{pp2.2}
Suppose that (A) holds.
Then for any $P,Q\in \mathcal{P}_2(\mathbb{R}^d)$ such that $\mathcal{S}(Q)$ is finite,
$V^S(\cdot,P,Q)\in C((0,T])$.
Besides, there exists a constant $\tilde{C}_T\ge 1$ which does not depend on $P,Q$
such that the following holds: for  $ t\in (0,T]$ and $\varepsilon\in (0,1)$,
\begin{equation}\label{2.4}
tV^S(t,P,Q)\le \tilde{C}_T \int_{\mathbb{R}^d}|x-y|^2P(dx)Q(dy)+
t\mathcal{S}(Q)+t\log (\tilde{C}_Tt^{d/2}),
\end{equation}
\begin{align}\label{2.21}
tV^S(t,P,Q)
&\ge \frac{1}{2\lambda_{\infty,t}}
\{(1-\varepsilon)^{2}{\bf T}_2(P,Q)\nonumber\\
&\qquad 
-\varepsilon^{-1}(1-\varepsilon)^{2}||\sigma ||_{\infty,t}^2t
-\varepsilon^{-1}(1-\varepsilon)|\xi|_{\infty,t}^2t^2\}.
\end{align}
In particular, $\{tV^S(t, P,Q)\}_{0< t\le T}$ is bounded, and 
\begin{align}
\limsup_{t\to 0}tV^S(t, P,Q)&\le \tilde{C}_T\int_{\mathbb{R}^d}|x-y|^2P(dx)Q(dy),\\
\liminf_{t\to 0}tV^S(t, P,Q)&\ge \frac{1}{2\lambda_{\infty,0+}}{\bf T}_2 (P,Q).
\end{align}
\end{proposition}

The following gives the long--time asymptotics of $\mu_t (P,Q)$ and $V^S(t,P,Q)$, and 
the upper bounds of $V^S(\cdot,P,Q)$ on $(0,\infty]$
when the solution to the SDE (\ref{1.4}) is ergodic.

\begin{theorem}\label{pp2.3}
Suppose that (A0) and (A4) hold.
Then for any $P,Q\in \mathcal{P}(\mathbb{R}^d)$, 
\begin{equation}\label{1.19L}
\lim_{t\to\infty}\mu_t (P,Q)=P\times Q,\quad \hbox{weakly.}
\end{equation}
Suppose, in addition, that $P,Q\in \mathcal{P}_2(\mathbb{R}^d)$ and $\mathcal{S}(Q)$ is finite.
Then
\begin{equation}\label{1.20}
\lim_{t\to\infty}V^S(t,P,Q)=H(Q(dy)|m(y)dy)<\infty.
\end{equation}
For any $T>0$, there exists a constant $\overline {C}_T >0$ 
which does not depend on $P,Q$ such that the following holds:
for $t\ge T$,
\begin{equation}\label{212.1}
V^S(t,P,Q)\le S(Q)+\overline {C}_T\left(1+\int_{\mathbb{R}^d}|x|^2(P(dx)+Q(dx))\right).
\end{equation}
\end{theorem}

The following holds from Theorem \ref{pp2.3}.

\begin{corollary}\label{co2.5}
Suppose that (A0) and (A4) hold.
Then for any  $P,Q\in \mathcal{P}_2(\mathbb{R}^d)$ such that  $\mathcal{S}(Q)$ is finite,
\begin{equation}\label{2.23}
\lim_{t\to\infty}H(P\times Q|\mu_t (P,Q))=0.
\end{equation}
\end{corollary}

\begin{remark}\label{Remark2.3}
(i) 
In \cite{Con},
they considered (\ref{1.19L})-(\ref{1.20}) on a smooth, connected, and complete Riemannian manifold without boundary.
In an Euclidean setting, their assumption is satisfied if $\xi(x)=-DU(x)$ for a convex function $U\in C^2 (\mathbb{R}^d)$
such that $D^2U$ is uniformly nondegenerate,
and if $H(Q|m)$ is finite and $P,Q$ have compact supports.
In particular, our assumption is weaker than theirs
though their estimates are much sharper than those in this paper.
They even gave the convergence rate in (\ref{1.20}).
(ii) (\ref{2.23}) implies the following by the CKP inequality (see \cite{Vi}) 
\begin{equation}
||\mu_t (P,Q)-P\times Q ||_{TV}\le \sqrt {2 H(P\times Q|\mu_t (P,Q))}\to 0,\quad t\to\infty,
\end{equation}
where $||\cdot ||_{TV}$ denotes the total variation distance (see e.g., \cite{Vi}).
\end{remark}

\section{Lemmas} \label{sec:3}

In this section, we give lemmas for the proof of our result.

We prove the following since we could not find appropriate literature.
\begin{lemma}\label{lemma3.1.1}
Let $Y,Z$ be $\mathbb{R}^d$-valued random variables defined on the same probability space.
Let $\{B(t)\}_{t\ge 0}$  be a $d$--dimensional standard Brownian motion that is independent of $(Y,Z)$.
Suppose that (A1,i) holds.
Then there exists a unique strong solution $\{X^{(Y,Z)}(t)\}_{0\le t< T}$ to the following:
\begin{equation}
X(t)=Y+\int_0^t \frac{Z-X(s)}{T-s}ds+ \int_0^t \sigma(s,X(s))dB(s), \quad 0\le  t<T.\label{3.1.0830}
\end{equation}
The following also holds.
\begin{equation}
\lim_{t\uparrow  T}X(t)=Z,\quad {\rm a.s.}\label{3.2.0830}.
\end{equation}
 In particular, $\{X^{(Y,Z)}(t)\}_{0\le t\le T}\in \mathcal{A}_T(P^Y,P^Z)$, where $X^{(Y,Z)}(T):=Z$.
\end{lemma}
\begin{proof} 
Under (A1,i), the  existence and the uniqueness of a strong solution to  (\ref{3.1.0830}) 
can be proved by a standard method (see e.g., \cite{IW14}).
We prove  (\ref{3.2.0830}).
The following holds: from (\ref{3.1.0830}),
\begin{equation}\label{3.3.0830}
X(t)-Z=(T-t)\left(\int_0^t\frac{\sigma(s,X(s))}{T-s}dB(s)+\frac{Y-Z}{T}\right),\quad 0\le t<T,
\end{equation}
since by It\^o's formula, 
$$d\frac{X(t)-Z}{T-t}=\frac{X(t)-Z}{(T-t)^2}dt+\frac{1}{T-t}\left(
\frac{Z-X(t)}{T-t}dt+\sigma(t,X(t))dB(t)\right).
$$
We prove the following to complete the proof:
\begin{equation}
\sup_{0\le t<T}\left\{(T-t)^{2/3}\left|\int_0^t\frac{\sigma(s,X(s))}{T-s}dB(s)\right|\right\}<\infty,\quad {\rm a.s.}.
\end{equation}
By It\^o's formula,
\begin{align*}
&(T-t)^{2/3}\int_0^t\frac{\sigma(s,X(s))}{T-s}dB(s)\\
&=-\frac{2}{3}\int_0^t(T-u)^{-1/3}\int_0^u \frac{\sigma(s,X(s))}{T-s}dB(s)du
+\int_0^t\frac{\sigma(s,X(s))}{(T-s)^{1/3}}dB(s),\\
\end{align*}
\begin{align*}
&E\left[\sup_{0\le t<T}\left|\int_0^t(T-u)^{-1/3}\int_0^u \frac{\sigma(s,X(s))}{T-s}dB(s)du\right|\right]\\
&\le E\left[\int_0^T(T-u)^{-1/3}\left|\int_0^u \frac{\sigma(s,X(s))}{T-s}dB(s)\right|du\right]\\
&\le \int_0^T(T-u)^{-1/3}E\left[\left|\int_0^u \frac{\sigma(s,X(s))}{T-s}dB(s)\right|^2\right]^{1/2}du\\
&=\int_0^T(T-u)^{-1/3}E\left[\int_0^u \frac{||\sigma(s,X(s))||^2}{(T-s)^2}ds\right]^{1/2}du\\
&\le ||\sigma||_{\infty,T}\int_0^T(T-u)^{-5/6}du=\frac{1}{6}||\sigma||_{\infty,T}T^{1/6}<\infty.
\end{align*}
By the monotone convergence theorem,
\begin{equation*}
E\left[\sup_{0\le t<T}\left|\int_0^t\frac{\sigma(s,X(s))}{(T-s)^{1/3}}dB(s)\right|^2\right]
=\lim_{u\uparrow T}E\left[\sup_{0\le t\le u}\left|\int_0^t\frac{\sigma(s,X(s))}{(T-s)^{1/3}}dB(s)\right|^2\right].
\end{equation*}
For $u\in (0,T)$, by  Doob's inequality,
\begin{align*}
&E\left[\sup_{0\le t\le u}\left|\int_0^t\frac{\sigma(s,X(s))}{(T-s)^{1/3}}dB(s)\right|^2\right]\\
&\le 4 E\left[\int_0^u\frac{||\sigma(s,X(s))||^2}{(T-s)^{2/3}}ds\right]\\
&\le 4||\sigma||_{\infty,T}^2\int_0^T(T-s)^{-2/3}ds=\frac{4}{3}||\sigma||_{\infty,T}^2T^{1/3}<\infty.
\end{align*}
\end{proof}

We give technical lemmas which will be made use of in the proofs of  Proposition  \ref{pp2.2} and Theorem \ref{pp2.3}.
We first describe three known results.

\begin{lemma}[see \cite{Csiszar,RT93}]\label{lm3.3}
Suppose that (A) holds.
Then for any $P,Q\in \mathcal{P}(\mathbb{R}^d), t\in (0,T]$ such that
$v^S(t,P,Q)$ is finite,
 $\mu_t(P,Q)$ is a unique minimizer of $v^S(t,P,Q)$
(see below (\ref{1.15}) for notation). 
\end{lemma}

\begin{lemma}[see \cite{A67}]\label{lm31}
Suppose that (A) holds.
Then for any $T>0$,
there exists $\tilde{C}_T\ge 1$ that is nondecreasing in T, and that
depends on $d, T, ||a||_{\infty,T}, |(D_x^*a)^*|_{\infty,T}, |\xi|_{\infty,T}$, and the infimum of the eigenvalues of $a$
such that the following holds: 
\begin{align}\label{31}
&\frac{1}{\tilde{C}_T(t-s)^{d/2}}\exp \left(-\frac{\tilde{C}_T|x-y|^2}{(t-s)}\right)\nonumber\\
&\le p(s,x;t,y)\nonumber\\
&\le \frac{\tilde{C}_T}{(t-s)^{d/2}}\exp \left(-\frac{|x-y|^2}{\tilde{C}_T(t-s)}\right),\quad x,y\in\mathbb{R}^d, 0\le s<t\le T.
\end{align}
\end{lemma}

\begin{lemma}[see Theorem 5 in \cite{V-2} and also \cite{V}\label{THM3.1}]
Suppose that (A0) and (A4) hold.
Then for any $k<(\ell-2)/2$,
\begin{equation}
\sup\left\{t^{k+1}(1+|x|^\ell)^{-1}\int_{\mathbb{R}^d} |p(t,x,y)-m(y)|dy:t> 0, x\in\mathbb{R}^d\right\}<\infty.
\end{equation}
\end{lemma}

Since
\begin{equation}\label{04021}
p(t+1,x,y)-m(y)=\int_{\mathbb{R}^d}(p(t ,x;z)-m(z))p(1,z,y)dz,\quad t>0,
\end{equation}
the following holds from Lemmas \ref{lm31}--\ref{THM3.1}.

\begin{lemma}\label{co31}
Suppose that (A0) and (A4) hold.
Then for any $k<(\ell-2)/2$,
\begin{equation}
\sup\{t^{k+1}(1+|x|^\ell)^{-1}|p(t +1,x,y)-m(y)|:t\ge 0, x,y\in\mathbb{R}^d\}<\infty.
\end{equation}
In particular, 
\begin{equation}
\lim_{t\to\infty}\left(\sup\{|p(t ,x,y)-m(y)|:y\in\mathbb{R}^d\}\right)=0, 
\end{equation}
locally uniformly in $x\in\mathbb{R}^d$.
\end{lemma}

We state and prove estimates on $p(t,x,y)$ and $m(y)$.

\begin{lemma}\label{LM3.3}
Suppose that (A0) holds.
Then for any $t\ge 0$ and $T>0$, the following holds: for $s,u\in (0,T]$
and $x,y\in\mathbb{R}^d$,
\begin{align}
p(s,x,y)
&\ge2^{-d/2}\tilde{C}_T^{-d-2}\exp \left(-\frac{2\tilde{C}_T|x|^2}{s}\right)p((2\tilde{C}_T^2)^{-1}s,0,y),\label{320}\\
p(s,x,y)
&\ge2^{-d/2}\tilde{C}_T^{-d-2}\exp\left(-\frac{2\tilde{C}_T|y|^2}{s}\right)p((2\tilde{C}_T^2)^{-1}s,x,0),
\label{321}
\end{align}
\begin{align}\label{317}
p(s+t+u,x,y)&\ge2^{-d}\tilde{C}_T^{-2(d+2)}\exp \left(-2\tilde{C}_T\left(\frac{|x|^2}{s}+\frac{|y|^2}{u}\right)\right)\nonumber
\\
&\qquad \times p((2\tilde{C}_T^2)^{-1}(s+u)+t,0,0).
\end{align}
\begin{equation}
p(t+T,x,y)\le\tilde{C}_TT^{-d/2}.\label{314}
\end{equation}
Suppose, in addition, that (A4) holds.
Then
\begin{align}\label{318}
\tilde{C}_TT^{-d/2}\ge m(y)&\ge2^{-d}\tilde{C}_T^{-2(d+2)}\exp \left(-\frac{2\tilde{C}_T|y|^2}{T}\right) m(0).
\end{align}
\end{lemma}

\begin{proof} 
We first prove (\ref{320})--(\ref{321}).
From Lemma \ref{lm31},
\begin{align}
&p(s,x,y)\nonumber\\
&\ge \frac{1}{\tilde{C}_Ts^{d/2}}\exp \left(-\frac{2\tilde{C}_T(|x|^2+|y|^2)}{s}\right)\nonumber\\
&=2^{-d/2}\tilde{C}_T^{-d-2}\exp \left(-\frac{2\tilde{C}_T|x|^2}{s}\right)
\frac{\tilde{C}_T}{\{(2\tilde{C}_T^2)^{-1}s\}^{d/2}}\exp \left(-\frac{|y|^2}{\tilde{C}_T\{(2\tilde{C}_T^2)^{-1}s\}}\right)\nonumber\\
&\ge2^{-d/2}\tilde{C}_T^{-d-2}\exp \left(-\frac{2\tilde{C}_T|x|^2}{s}\right)p((2\tilde{C}_T^2)^{-1}s,0,y).
\end{align}
(\ref{321}) can be proved in the same way.
(\ref{320})--(\ref{321}) imply (\ref{317}) by Chapman--Kolmogorov's equality: for $t>0$,
\begin{align}
p(s+t+u,x,y)&=\int_{\mathbb{R}^d\times \mathbb{R}^d}p(s,x,z_1)p(t,z_1,z_2)p(u,z_2,y)dz_1dz_2,\\
p(s+u,x,y)&=\int_{\mathbb{R}^d}p(s,x,z)p(u,z,y)dz.\label{319}
\end{align}
Replace $s$ and $u$ by $t$ and $T$  in  (\ref{319}), respectively.
Lemmas \ref{lm31} implies (\ref{314}).
Let $x=0, u=T$ in  (\ref{317}) and then $t\to\infty$ in (\ref{317})--(\ref{314}).
Then we obtain (\ref{318}) from Lemma \ref{co31}.
\end{proof}

Under (A0), for $P, Q\in\mathcal{P}(\mathbb{R}^d)$ and $T>0$,
\begin{align}\label{323.}
\varphi_1(t,y)&:=\log\left(\int_{\mathbb{R}^d}p(t,x,y)\nu_1(dx)\right),\quad (t,x)\in (0,T]\times \mathbb{R}^d,\\
\varphi_2(t,x)&:=\log\left(\int_{\mathbb{R}^d}p(T-t,x,y)\nu_2(dy)\right),\quad (t,x)\in [0,T)\times \mathbb{R}^d
\end{align}
(see (\ref{113}) for notation).
Then $\varphi_1\in C^{1,2}((0,T]\times \mathbb{R}^d), \varphi_2\in C^{1,2}([0,T)\times \mathbb{R}^d)$ (see \cite{J75}).
By definition of $\mu_T=\mu_T(P,Q)\in \Pi (P,Q)$, 
\begin{equation}\label{324.}
\mu_T(dxdy)=\exp(-\varphi_1(T,y)-\varphi_2(0,x))p(T,x,y)P(dx)Q(dy).
\end{equation}

Since we could not find appropriate literature,
we prove the following for the sake of completeness.
\begin{lemma}\label{LM3.9}
Suppose that (A0) holds.
Then for any $P, Q\in\mathcal{P}_2(\mathbb{R}^d)$ such that 
$\mathcal{S}(Q)$ is finite, $ \varphi_1(T,\cdot)\in  L^1 (Q)$ and $\varphi_2(0,\cdot)\in L^1 (P)$.
\end{lemma}

\begin{proof}
From (\ref{320})--(\ref{321}),  
\begin{equation}\label{3.28}
\min (\varphi_1(T,x),\varphi_2(0,x))\ge-C_1|x|^2-C_2,\quad x\in \mathbb{R}^d,
\end{equation}
where
\begin{align*}
C_1&:=\frac{2\tilde{C}_T}{T}>0,\\
C_2&:=\log (2^{d/2}\tilde{C}_T^{d+2})\\
&\qquad -\min (0, \varphi_1((2\tilde{C}_T^2)^{-1}T,0),\varphi_2((1-(2\tilde{C}_T^2)^{-1})T, 0))>0.
\end{align*}
Recall that $\tilde{C}_T\ge 1$ (see Lemmas \ref{lm31}).
From (\ref{3.28}),
\begin{align}
&\int_{\mathbb{R}^d} |\varphi_1(T,y)|Q(dy)+\int_{\mathbb{R}^d} |\varphi_2(0,x)|P(dx)\nonumber
\\
&\le \int_{\mathbb{R}^d} (\varphi_1(T,y)+2(C_1|y|^2+C_2))Q(dy)\nonumber\\
&\qquad+\int_{\mathbb{R}^d} (\varphi_2(0,x)+2(C_1|x|^2+C_2))P(dx)\nonumber\\
&=-V^S(T,P,Q)+\mathcal{S}(Q)+
2\int_{\mathbb{R}^d\times \mathbb{R}^d} (C_1(|x|^2+|y|^2)+2C_2)\mu_T(dxdy)<\infty
\end{align}
since $\mathcal{S}(Q)<\infty$ and $P, Q\in\mathcal{P}_2(\mathbb{R}^d)$.
\end{proof}


\section{Proof}\label{sec:4}
In this section, we prove our result.

\begin{proof} [Proof of  Proposition~{\upshape\ref{Thm2.1}}]

We only have to consider the case where ${\bf T}_r(P,Q)$ is finite.
Replace $T$ by $t\in (0,T]$ in Lemma \ref{lemma3.1.1}
and take $Y,Z$ such that $E[|Y-Z|^r]$ is finite and that $\{X^{(Y,Z)}(s)\}_{0\le s\le t}\in \mathcal{A}_t(P,Q)$.
Then from (\ref{3.3.0830}),
\begin{align}
&V_r(t,P,Q)\nonumber\\
&\le E\left[\int_0^t \left|\frac{Z-X^{(Y,Z)}(s)}{t-s}\right|^rds\right]\nonumber\\
&\le \int_0^t E\left[E\left[\left|\frac{Z-X^{(Y,Z)}(s)}{t-s}\right|^2\biggl|Y,Z\right]^{r/2}\right]ds\nonumber\\
&=\int_0^t E\left[E\left[\int_0^s\frac{||\sigma(u,X^{(Y,Z)}(u))||^2}{(t-u)^2}du+\frac{|Y-Z|^2}{t^2}\biggl|Y,Z\right]^{r/2}\right]ds\nonumber\\
&\le\frac{2||\sigma||_{\infty,t}^r}{2-r}t^{1-r/2}+E[|Y-Z|^r] t^{1-r}.
\end{align}
Indeed,
\begin{align*}
&E\left[\int_0^s\frac{||\sigma(u,X^{(Y,Z)}(u))||^2}{(t-u)^2}du+\frac{|Y-Z|^2}{t^2}\biggl|Y,Z\right]^{r/2}\\
&\le E\left[\int_0^s\frac{||\sigma||_{\infty,t}^2}{(t-u)^2}du+\frac{|Y-Z|^2}{t^2}\biggl|Y,Z\right]^{r/2}\\
&\le \frac{||\sigma||_{\infty,t}^r}{(t-s)^{r/2}}+\frac{|Y-Z|^r}{t^r}
\end{align*}
since
\begin{equation}\label{4.9.m1}
(a+b)^{r}\le a^{r}+b^{r}, \quad a,b\ge 0, 0\le r\le 1.
\end{equation}
\end{proof}

We prove Proposition \ref{pp2.2.0}.

\begin{proof} [Proof of Proposition~{\upshape\ref{pp2.2.0}}]
We only have to consider the case where $V_r(t, P,Q)$ is finite.
Take $X\in\mathcal{A}_t( P,Q)$ such that 
$$E\left[ \int_0^t |u_X(s)|^r ds\right]<\infty.$$
(\ref{2.7.0})--(\ref{2.6}) can be proved by the following.
By Jensen's inequality, 
\begin{align}
t^{r-1}E\left[ \int_0^t |u_X(s)|^r ds\right]
&\ge E\left[\left|\int_0^t u_X(s)ds\right|^r\right]\nonumber\\
&=E\left[\left|X (t)- X(0)-\int_0^t \sigma (s, X(s))dB(s)\right|^r\right]
\end{align}
By Minkowski's inequality, the following also holds:
\begin{align}
&E\left[\left|X (t)- X(0)-\int_0^t \sigma (s, X(s))dB(s)\right|^r\right]^{1/r}\nonumber\\
&\ge E\left[\left|X (t)- X(0)\right|^r\right]^{1/r} -E\left[\left|\int_0^t\sigma (s, X(s))dB(s)\right|^r\right]^{1/r}.
\end{align}
By Schwartz's and Burkholder--Davis--Gundy's inequalities, there exists $C_r>0$ such that the following holds:
\begin{align}
&E\left[\left|\int_0^t \sigma (s, X(s))dB(s)\right|^r\right]\nonumber\\
&\le 
\begin{cases}
\displaystyle E\left[\int_0^t ||\sigma (s, X(s))||^2ds\right]^{r/2},&
1\le r\le 2,\\
\displaystyle 
C_rE\left[\left(\int_0^t ||\sigma (s, X(s))||^2ds\right)^{r/2}\right],&
r\ge 1,
\end{cases}
\end{align}
which implies  (\ref{2.7.0}). (\ref{2.7.0}) implies (\ref{2.6})
since
\begin{align}\label{4.6.0905}
|a+b|^r&=\left|(1-\varepsilon)\frac{a}{1-\varepsilon}+\varepsilon \frac{b}{\varepsilon}\right|^r\nonumber\\
&\le (1-\varepsilon)^{1-r}|a|^r+\varepsilon^{1-r}|b|^r,\quad a,b\in \mathbb{R}^d, r\ge 1.
\end{align}
Without loss of generality, we can assume that $X^0(0)=X(0)$.
Indeed, $X^0$ is a strong solution of $(\ref{1.9})$ with $b_X=0$ for any starting point
from (A1,i) (see e.g., \cite{IW14}).
(\ref{2.5})--(\ref{2.4.0}) can be proved in the same way as above since
\begin{equation}
\int_0^t u_X(s)ds=X (t)- X^0(t)+\int_0^t (\sigma (s, X^0(s))-\sigma (s, X(s)))dB(s).
\end{equation}
Letting $\varepsilon\to 0$ in (\ref{2.4.0}), we obtain 
(\ref{2.6.0}).
\end{proof}

We prove Corollary \ref{co2.3}.
\begin{proof} [Proof of Corollary~{\upshape\ref{co2.3}}]
Since $r\in [1,2)$, we use (\ref{2.7.0})--(\ref{2.6}) with $C_r$ replaced by $1$.
We first prove (\ref{2.7}).
In (\ref{2.6}), set $\varepsilon=Rt^{1/2}, R<t^{-1/2}$.
Then 
\begin{align}\label{4.4.0830}
&t^{-1/2}(t^{r-1}V_r(t,P,Q)-{\bf T}_r(P,Q))\nonumber\\
&\ge -R\frac{(1-\varepsilon)^{r-1}-1}{-\varepsilon}{\bf T}_r(P,Q)
-R^{1-r}(1-\varepsilon)^{r-1}||\sigma ||_{\infty,t}^r\nonumber\\
&\to -R(r-1){\bf T}_r(P,Q)-R^{1-r}||\sigma ||_{\infty,0+}^r,\quad t\to 0.
\end{align}
When $r=1$, the proof is over from (\ref{2.2.0527}).
When $r>1$, taking the maximum of the r.h.s. of (\ref{4.4.0830}) in $R>0$, we obtain  (\ref{2.7})
from (\ref{2.2.0527}).
Indeed, for $c>0$,
$$x+cx^{1-r}\ge (c(r-1))^{1/r}+c(c(r-1))^{(1-r)/r},\quad x>0,$$
where the equality holds if and only if $x= (c(r-1))^{1/r}$.
If ${\bf T}_r(P,Q)>0$, then
\begin{align*}
&R(r-1){\bf T}_r(P,Q)+R^{1-r}||\sigma ||_{\infty,0+}^r\\
&=(r-1){\bf T}_r(P,Q)\left \{R+\frac{||\sigma ||_{\infty,0+}^r}{(r-1){\bf T}_r(P,Q)}R^{1-r}\right\}\\
&\ge (r-1){\bf T}_r(P,Q)\left \{\left(\frac{||\sigma ||_{\infty,0+}^r}{{\bf T}_r(P,Q)}\right)^{1/r}+\frac{||\sigma ||_{\infty,0+}^r}{(r-1){\bf T}_r(P,Q)}\left(\frac{||\sigma ||_{\infty,0+}^r}{{\bf T}_r(P,Q)}\right)^{(1-r)/r}\right\}\\
&=r||\sigma ||_{\infty,0+}{\bf T}_r(P,Q)^{1-1/r}.
\end{align*}
If ${\bf T}_r(P,Q)=0$, then
$$
R(r-1){\bf T}_r(P,Q)+R^{1-r}||\sigma ||_{\infty,0+}^r\to 0, \quad R\to\infty.
$$
(\ref{2.8}) can be proved from (\ref{2.2.0527}).
We prove (\ref{2.11}).
From (\ref{2.2.0527}) and (\ref{4.9.m1}),
\begin{align}
(t^{r-1}V_r(t, P,Q))^{1/r}&\le \left({\bf T}_r (P,Q)+\frac{2||\sigma||_{\infty,t}^r }{2-r}t^{r/2}\right)^{1/r}\nonumber\\
&\le{\bf T}_r (P,Q)^{1/r}+\left(\frac{2}{2-r}\right)^{1/r} ||\sigma||_{\infty,t} t^{1/2}.
\end{align}
(\ref{2.7.0}) together with the following completes the proof of  (\ref{2.11}).
$$\frac{2}{2-r}> 1,\quad 1\le r<2.$$
\end{proof}

We prove Theorem \ref{Theorem2.1}.

\begin{proof} [Proof of Theorem~{\upshape\ref{Theorem2.1}}]

We only have to consider the case where $V_r(t,P,Q)$ is finite.
Take $X\in\mathcal{A}_t( P,Q)$ such that 
$$E\left[ \int_0^t |u_X(s)|^r ds\right]<\infty.$$
By Jensen's inequality, 
\begin{equation}
E\left[ \int_0^t |u_X(s)|^r ds\right]
\ge t^{1-r} E\left[\left|\int_0^t u_X(s)ds\right|^r\right].
\end{equation}
Since $u\mapsto |u|^r$ is convex,
\begin{align}
\left|\int_0^t u_X(s)ds\right|^r
&=\left|X (t)- X(0)-\int_0^t \sigma (s, X(s))dB(s)\right|^r\nonumber\\
&\ge \left|-\int_0^t \sigma (s, X(s))dB(s)\right|^r\nonumber\\
&\qquad -r\left|-\int_0^t \sigma (s, X(s))dB(s)\right|^{r-1}|X (t)- X(0)|.
\end{align}
\begin{align}
&E\left[\left|-\int_0^t \sigma (s, X(s))dB(s)\right|^{r-1}|X (t)- X(0)|\right]\nonumber\\
&\le E\left[\left|-\int_0^t \sigma (s, X(s))dB(s)\right|^r\right]^{(r-1)/r}
E[|X (t)- X(0)|^r]^{1/r}\nonumber\\
&\le \{(\lambda_M t)^{r/2}\}^{(r-1)/r}\times 
\left(2^{r-1}\int_{\mathbb{R}^d}|x|^r(P(dx)+Q(dx))\right)^{1/r}.
\end{align}
By Schwartz's inequality,
\begin{align}
&E\left[ \left|-\int_0^t \sigma (s, X(s))dB(s)\right|^{r}\right]\nonumber\\
&\ge E\left[ \left|-\int_0^t \sigma (s, X(s))dB(s)\right|^{3r/2}\right]^2\nonumber\\
&\qquad \times 
E\left[ \sup_{0\le s\le t}\left|-\int_0^s \sigma (u, X(u))dB(u)\right|^{2r}\right]^{-1}
\end{align}
since 
\begin{align*}
&\left|-\int_0^t \sigma (s, X(s))dB(s)\right|^{3r/2}\\
&\le \left|-\int_0^t \sigma (s, X(s))dB(s)\right|^{r/2}\times 
\sup_{0\le s\le t}\left|-\int_0^s \sigma (u, X(u))dB(u)\right|^{r}.
\end{align*}
The following completes the proof:
by Doob's and Burkholder--Davis--Gundy's inequalities, there exists $c_{3r/4}, C_r>0$ such that the following holds:
\begin{align*}
&E\left[ \left|-\int_0^t \sigma (s, X(s))dB(s)\right|^{3r/2}\right]^2\nonumber\\
&\ge \left(\frac{3r/2}{3r/2-1}\right)^{-3r}
E\left[ \sup_{0\le s\le t}\left|-\int_0^s \sigma (u, X(u))dB(u)\right|^{3r/2}\right]^2\nonumber\\
&\ge \left(\frac{3r/2}{3r/2-1}\right)^{-3r}
c_{3r/4}^{2}
E\left[\left(\int_0^t ||\sigma (u, X(u))||^2ds\right)^{3r/4}\right]^2\nonumber\\
&\ge \left(\frac{3r/2}{3r/2-1}\right)^{-3r}c_{3r/4}^{2}(\lambda_mt)^{3r/2},
\end{align*}
\begin{align*}
E\left[ \sup_{0\le s\le t}\left|-\int_0^s \sigma (u, X(u))dB(u)\right|^{2r}\right]
&\le C_rE\left[\left(\int_0^t ||\sigma (u, X(u))||^2ds\right)^{r}\right]\nonumber\\
&\le C_r (\lambda_Mt)^{r},
\end{align*}
\end{proof}

Next, we prove Proposition \ref{pp2.1}.

\begin{proof}[Proof of  Proposition~{\upshape\ref{pp2.1}}]
We prove the first part of the proposition.
Let $Y,Z$ be $\mathbb{R}^d$--valued random variables defined on the same probability space
such that $P^{(Y,Z)}\in \Pi (P,Q)$ and that $E[|Y-Z|^r]$ is finite.
Let $\{B(t)\}_{t\ge 0}$  be a $d$--dimensional standard Brownian motion that is independent of $(Y,Z)$.
Construct $\{X^0(t)\}_{0\le t\le T}$ so that $X^0 (0)=Y$.
Under (A1,i), for $\delta\in (0,T)$, let $\{X^\delta (t)\}_{0\le t\le T}\in \mathcal{A}_T(P,Q)$  be a solution to the following (see Lemma \ref{lemma3.1.1}):
\begin{align} 
X^\delta (t)&:=X^0 (t),\quad  0\le t\le T-\delta,\nonumber\\
dX^\delta (t)&=
\frac{Z-X^\delta (t)}{T-t}dt+\sigma (t, X^\delta (t))dB(t),\quad T-\delta < t<T.
\end{align}
Notice that $\{B(t)-B(T-\delta)\}_{t\ge T-\delta}$ and $\{Z, X^0 (t),0\le t\le T-\delta\}$ are mutually independent.
From (A3),
\begin{align}\label{4.11.0905} 
V(T,P,Q)&\le E\biggl[\int_{T-\delta}^T L\left(t,X^\delta (t) ;\frac{Z-X^\delta (t)}{T-t}\right)dt \biggr]\nonumber\\
&\le C_{r,T}E\biggl[\int_{T-\delta}^T \left|\frac{Z-X^\delta (t)}{T-t}\right|^rdt \biggr]
+\delta C_{r,T}'\nonumber\\
&\to 0,\quad \delta\to 0.
\end{align}
Indeed, in  the same way as in the proof of Proposition \ref{Thm2.1}, 
$$
E\biggl[\int_{T-\delta}^T \left|\frac{X^\delta (t)-Z}{T-t}\right|^rdt \biggr]
\le \frac{2||\sigma||_{\infty,T}^r }{2-r}\delta^{1-r/2}+
\delta^{1-r}E[|X^0 (T-\delta)-Z|^r],
$$
$$|X^0 (T-\delta)-Z|^r \le |X^0 (T-\delta)-X^0 (0)|^r+|Y-Z|^r,$$
from (\ref{4.9.m1}).
By Schwartz's inequality,
$$E[|X^0 (T-\delta)-X^0 (0)|^r]
\le E\left[\left(\int_0^{T-\delta} ||\sigma (s, X^0(s))||^2ds\right)\right]^{r/2}
\le ||\sigma||_{\infty,T}^r T^{r/2}.
$$

We prove the second part of the proposition.
Since $V_r(T, P,Q)$ is finite from the assumption, there exists 
$\{X(t)\}_{0\le t\le T}\in \mathcal{A}_T( P,Q)$ such that 
$$E\biggl[\int_0^T |u_X(t)|^rdt\biggr]<\infty.$$
For $n\ge T^{-1}$, 
\begin{align} 
u_n(t)&:=
\begin{cases}
0,&0\le t\le T-\frac{1}{n},\\
nT\times u_X(nT(t-T+\frac{1}{n})), &T-\frac{1}{n}< t\le T,
\end{cases}\nonumber\\
X_n (t)&:=X(0)+\int_0^t u_n(s)ds+\int_0^t\sigma (s)dB(s),\quad 0\le t\le T.
\end{align}
Then 
\begin{align}
(X_n (0),X_n (T))&=(X(0),X(T)),\nonumber\\
nT\left(t-T+\frac{1}{n}\right)&\le t,\quad 0\le t\le T.
\end{align}
In particular, $X_n\in\mathcal{A}_T(P,Q)$ and the following holds: 
\begin{equation}
V(T,P,Q)\le E\biggl[\int_0^T L(t,X_n (t);u_n(t))dt \biggr]\to 0, \quad n\to\infty.
\end{equation}
Indeed, from (A3), 
\begin{align*}
E\biggl[\int_0^T L(t,X_n (t);u_n(t))dt \biggr]
&\le E\biggl[C_{r,T}\int_{T-1/n}^T |u_n(t)|^rdt+C_{r,T}'n^{-1}\biggr]\\
&=  C_{r,T}(nT)^{r-1}E\biggl[\int_{0}^T |u_X(t)|^rdt\biggr]+C_{r,T}'n^{-1}.
\end{align*}
\end{proof}

We prove Proposition \ref{pp2.2}.

\begin{proof} [Proof of Proposition~{\upshape\ref{pp2.2}}]
Since $V^S(t,P,Q)=v^S(t,P,Q)$ under (A)  (see section \ref{intro}), we consider $v^S$.
Under (A), $\mathbb{R}^d\times (0,T]\times \mathbb{R}^d\ni (x,t,y)\mapsto p(0,x;t,y)$ is 
positive and continuous (see Theorem \ref{thm11}).
In particular, $ (0,T]\ni t\mapsto p(0,x;t,y)$ is continuous locally uniformly in $(x,y)$.
From \cite{mikami2021}, $ (0,T]\ni t\mapsto \mu_t (P,Q)$ is weakly continuous.
Since the relative entropy $H(\mu|\nu)$ is lower semicontinuous in $\mu,\nu$ (see e.g., \cite{DE}), 
the following holds: 
\begin{align}\label{322}
&\liminf_{s\to t }H(\mu_s (P,Q)|P(dx)p(0,x;s,y)dy)\nonumber\\
&\ge H(\mu_t (P,Q)|P(dx)p(0,x;t,y)dy)=v^S(t,P,Q), \quad \hbox {on }(0,T].
 \end{align}
 Since $P,Q\in \mathcal{P}_2(\mathbb{R}^d)$,
from Lemma \ref{lm31}, by the dominated convergence theorem, 
\begin{equation}\label{4.17.J}
\limsup_{s\to t }v^S(s,P,Q)\le v^S(t,P,Q), \quad \hbox {on }(0,T].
 \end{equation}
Indeed, from (\ref{1.15}) and Lemma \ref{lm3.3},
\begin{align}\label{323}
v^S(s,P,Q)
&=H(\mu_s (P,Q)|P(dx)p(0,x;s,y)dy)\nonumber\\
&\le H(\mu_t (P,Q)|P(dx)p(0,x;s,y)dy)\nonumber\\
&=H(\mu_t(P,Q)|P(dx)p(0,x;t,y)dy)\nonumber\\
&\qquad 
+\int_{\mathbb{R}^d\times \mathbb{R}^d}
\left(\log\frac{p(0,x;t,y)}{p(0,x;s,y)} \right)\mu_t(P,Q)(dxdy). 
\end{align}
Theorem \ref{thm1.1} implies that $\{tv^S(t, P,Q)\}_{0<t\le T}$ is bounded and hence $v^S(\cdot, P,Q)\in C((0,T])$. 
(\ref{2.4}) can be proved from Lemma \ref{lm31}: in the same way as (\ref{323}),
\begin{align}\label{324}
v^S(t,P,Q)&\le H(P\times Q|P(dx)p(0,x;t,y)dy)\nonumber\\
&=\mathcal{S}(Q) -\int_{\mathbb{R}^d\times \mathbb{R}^d}
\left(\log p(0,x;t,y)\right)P(dx)Q(dy). 
\end{align}
(\ref{2.6}) together with the following implies (\ref{2.21}):
for $u,x\in \mathbb{R}^d,s\in [0,t], \varepsilon \in (0,1)$,
\begin{equation}
|u|^2
\le (1-\varepsilon )^{-1}\lambda_{\infty,t}|\sigma(s,x)^{-1}(u-\xi(s,x))|^2
+ \varepsilon^{-1} |\xi|_{\infty,t}^2.
\end{equation}
Indeed,  for $\xi\in \mathbb{R}^d$,
$$|u|^2=|(1-\varepsilon )(1-\varepsilon )^{-1}(u-\xi)+\varepsilon  \varepsilon ^{-1}\xi|^2
\le (1-\varepsilon )^{-1}|u-\xi|^2+ \varepsilon^{-1}|\xi|^2,$$
$$|u-\xi|^2= |\sigma(s,x)\sigma(s,x)^{-1}(u-\xi)|^2
\le \lambda_{\infty,t}|\sigma(s,x)^{-1}(u-\xi)|^2.$$
\end{proof}

We prove Theorem \ref{pp2.3}.

\begin{proof}  [Proof of Theorem~{\upshape\ref{pp2.3}}]
We prove (\ref{1.19L}).
Under  (A0) and (A4), $p(t,x,y)$ and $m(y)$ are positive and continuous (see Theorem \ref{thm11}, Remark \ref{remark21} and Lemma \ref{lm31}).
From Lemma \ref{co31} and \cite{mikami2021}, 
there exists a unique nonnegative $\sigma$--finite measure $\tilde\nu_1(dx)\tilde\nu_2(dy)$ such that 
the following holds:
\begin{align}\label{3.15.m}
\lim_{t\to\infty}\mu_t (P,Q)&=\tilde\nu_1(dx)m(y)\tilde\nu_2(dy),\quad \hbox{weakly},\nonumber\\
\tilde\nu_1(dx)m(y)\tilde\nu_2(dy)&\in \Pi (P,Q).
\end{align}
Since $\tilde\nu_1(dx)m(y)\tilde\nu_2(dy)\in \Pi (P,Q)$, 
\begin{equation}
\tilde\nu_1(dx)m(y)\tilde\nu_2(dy)=P(dx)Q(dy).
 \end{equation}

We prove (\ref{1.20}) by considering $v^S(t,P,Q)$.
Since $Q\in \mathcal{P}_2(\mathbb{R}^d)$ and  $\mathcal{S}(Q)$ is finite,
Lemmas \ref{LM3.3} implies that $ H(Q(dy)|m(y)dy)$ is finite.
In the same way as (\ref{322}), from (\ref{1.19L}) and Lemma \ref{co31},
\begin{equation}\label{3.24}
\liminf_{t\to\infty}H(\mu_t (P,Q)|P(dx)p(t,x,y)dy)\ge H(Q(dy)|m(y)dy).
 \end{equation}
 We prove the following:
 \begin{equation}\label{421}
\limsup_{t\to\infty}H(\mu_t (P,Q)|P(dx)p(t,x,y)dy)\le H(Q(dy)|m(y)dy).
 \end{equation}
From (\ref{324}),
\begin{equation}\label{3.25}
v^S(t,P,Q)\le H(Q(dy)|m(y)dy)-\int_{\mathbb{R}^d\times \mathbb{R}^d}
\left(\log\frac{p(t,x,y)}{m(y)} \right)P(dx)Q(dy).
\end{equation}
From Lemma \ref{LM3.3}, for $T,t>0$, 
\begin{align}\label{4.17}
&p(t+2T,x,y)m(y)^{-1}\nonumber\\
&\ge 2^{-d}\tilde{C}_T^{-2(d+2)}\exp \left(-2\tilde{C}_T\left(\frac{|x|^2}{T}+\frac{|y|^2}{T}\right)\right)
p((2\tilde{C}_T^2)^{-1}\times 2T+t,0,0)\nonumber\\
&\times \tilde{C}_T^{-1}T^{d/2}.
\end{align}
Letting $t\to\infty$, Lemma \ref{co31} and Fatou's lemma imply  (\ref{421}).


We prove (\ref{212.1}).
From (\ref{324}), 
the following completes the proof.
Setting $s=u=T/2$ in   (\ref{317}),
the following holds: for any $t\ge 0$,
\begin{align}\label{328}
&\log p(t+T,x,y) \nonumber\\
&\ge\log\left(2^{-d}\tilde{C}_T^{-2(d+2)} 
p((2\tilde{C}_T^2)^{-1}T+t,0,0)\right)-\frac{4\tilde{C}_T}{T}(|x|^2+|y|^2).
\end{align}
From Lemmas \ref{lm31} and \ref{co31}, 
\begin{equation}
p((2\tilde{C}_T^2)^{-1}T+t,0,0)\ge \inf\{p(u,0,0):u\ge (2\tilde{C}_T^2)^{-1}T\}>0,\quad t\ge 0.
\end{equation}
Indeed, $\{ p(u,0,0):u\ge (2\tilde{C}_T^2)^{-1}T\}$ is locally uniformly positive, $p(u,0,0)\to m(0)$ as $u\to\infty$ and $m(0)>0$.
\end{proof}

We prove Corollary \ref{co2.5}.

\begin{proof} [Proof of  Corollary~{\upshape\ref{co2.5}}]
For $t>0$,
\begin{align}\label{4.25}
\quad 0&\le H(P\times Q|\mu_t(P, Q))\nonumber\\
&=-V^S(t,P,Q)+H(Q(dy)|m(y)dy)\nonumber\\
&\qquad -\int_{\mathbb{R}^d\times \mathbb{R}^d} \left(\log \frac{p(t,x,y)}{m(y)}\right)P(dx)Q(dy).
\end{align}
Indeed, from (\ref{324.}) and Lemma \ref{LM3.9},
\begin{align*}
&H(P\times Q|\mu_t)\\
&=\int_{\mathbb{R}^d\times \mathbb{R}^d} (\varphi_1(t,y)+\varphi_2(0,x)-\log p(t,x,y))P(dx)Q(dy)\nonumber\\
&=-\int_{\mathbb{R}^d\times \mathbb{R}^d} (\log q(y)-\varphi_1(t,y)-\varphi_2(0,x))\mu_t(dxdy)\nonumber\\
&\qquad +\int_{\mathbb{R}^d}\left(\log \frac{q(y)}{m(y)} \right)q(y)dy
-\int_{\mathbb{R}^d\times \mathbb{R}^d} \left(\log \frac{p(t,x,y)}{m(y)}\right)P(dx)Q(dy)\nonumber
\end{align*}
(see Lemmas \ref{lm31} and \ref{LM3.3} for the upper and lower bounds of $p(t,x,y)$ and $m(y)$).

The limsup of the r.h.s. of (\ref{4.25}) as $t\to\infty$ is less than or equal to $0$, from (\ref{1.20}), Lemma \ref{co31} and  (\ref{4.17}) by Fatou's lemma. 
\end{proof}


\end{document}